\documentclass[10pt,twoside]{article}
\usepackage{Latex-document}

\title{\bf Reforms of the University Mathematics \vskip -2mm
Education for Non-mathematical \vskip -2mm Specialties \vskip 6mm}

\author{ Shutie Xiao\vspace*{-0.5cm}\thanks{Department of Applied
Mathematics, Tsinghua University, South Building 11--5--402,
Beijing 100084, China. E-mail: xstwef@mail.tsinghua.edu.cn}}

\date{\vspace{-8mm}}

\begin{document}

\maketitle

\thispagestyle{first} \setcounter{page}{897}

\begin{abstract}

\vskip 3mm

This article is a part of the report for the research project ``Reform of the Course System and Teaching Content
of Higher Mathematics (For Non-Mathematical Specialties)'' in 1995, supported by the National Ministry of
Education. There are thirteen universities participated in this project. The Report not only reflects results of
our participants, but also includes valuable opinions of many colleagues in the mathematical education circles.

In this article, after a brief description on the history and reform situation of the higher mathematics education
in China, attention concentrates to three aspects. They are: main problems in this field existing; the functions
of mathematics accomplishment for college students; these concern course system, teaching and learning philosophy,
such as overemphasized specialty education, overlooking to arouse rational thinking and aesthetic conceptions,
etc. The last aspect contains a discussion on several important relationships, such as: knowledge impartment and
quality cultivation, inherence and modernization of the mathematical knowledge, teacher's guidance and students'
initiative, mathematical basic training and mathematical application consciousness and ability cultivation, etc.

\vskip 4.5mm


\end{abstract}

\vskip 12mm

\section{Introduction}

\vskip-5mm \hspace{5mm}

This is a part of the report for the research project of `` Reforms on the Course System and Teaching Content of
Higher Mathematics (For Non-Mathematical Specialties)'' in 1995, supported by National Ministry of Education.
There are thirteen universities participated in this research project. They worked together cooperatively.
Besides, many experts and scholars in mathematical and educational circles had been consulted in many aspects.
This essay not only reflects the research results of our participants-universities in more than four years, but
also includes valuable opinions of many colleagues in the mathematical education circle.

\section{Reform situation of the university mathematics education
in China and the main problems existing}

\vskip-5mm \hspace{5mm}

For the past half a century, our modern university education, accompanying the progress process of our socialist
political and economic development, experiencing the stage of ``learning from the Soviet Union in all respects'',
many times of ``educational reforms'' and today's ``reform and opening to the world'', has grown in wave-like
development and formed today's scale. In this process, the teaching reform on the university non-mathematics class
specialties higher mathematics has almost not been stopped. It can be said that mathematics teaching in this
period, generally speaking, adapted itself to the requirements for training professionals in various fields under
the planned economic system and made important contributions for the society. In this respect, several generations
of educators of mathematical teaching devoted their wisdom and hardworking.

Historically, the focus of the reform of the university mathematical education in China is always concentrated on
the problem of the combination of mathematics with practice. In a long period before the reforms and opening to
the world, the basic task of the non-mathematics class specialties mathematical teaching is for specialties, or
mathematics courses should serve courses of specialties. As mentioned earlier, in the tendency of the extreme
situation of putting too much emphasis on specialties, such a saying is natural; but we went even further. In
addition, we were in a situation of long-period of isolation. We didn't know the situation of mathematics
education in other countries, even were not clear on the reform in the Soviet Union. Without comparison, it is
difficult to discover the problems. Moreover, different opinions on such a problem as ``the relationship between
the teaching of fundamental mathematics and that of specialties,'' should have been a problem of academic views in
the first place and can be discussed and tested; but such opinions were always improperly treated as ideological
problems even the political ones. Therefore, it was difficult to express different views. Since the 1980s, our
policies of the reform and opening to the world have brought spring for the educational reforms. When we learned
of the changes in those decades of years, we realized that we had lagged behind and we should do our best to catch
up. At that time, people found that, as a developing country, China was confronted with a series of serious
challenges in the socialist market economy. In addition, the rapid appearance of the knowledge economy based on
the intelligent resources and innovative competition, we are forced to reconsider problems, such as educational
concepts, teaching systems, course systems, teaching content and methods etc., which are formed in the original
planned economy.

Viewing from the angle of the university mathematics education,
the existing main problems may be summarized as follows:

First, the ``specialty education'' is overemphasized and thus a one-sided understanding of the role of the
university mathematics education ``serving specialties'' is formed. Such an understanding, as the guiding idea of
education, is embodied in every link of teaching. Due to its long duration, the depth and extension of its effects
cannot be underestimated. Even now, among quite a large part of teaching cadres and teachers, such an
understanding may still work, because in the past the problem of ``what is the role of university mathematical
courses in the university education'' had not been studied properly.

Secondly, as the cultivation target was training ``specialists and engineers who could work directly after their
graduation,'' the teaching process was speeded up. Moreover, we have the tradition of putting emphasis on the
classroom teaching, and the course of mathematics itself has its own specific strict logical system, the
instillation-type of teaching methods can be said to be deep-rooted.

Thirdly, as mathematics teachers of non-mathematical class specialties
usually undertake heavy teaching tasks, and in the technological and
engineering colleges and universities, the mathematical scientific research
not being combined with practice directly is viewed
as being isolated from practice, the teachers engaged in
teaching basic courses of mathematics have not touched scientific research
for a long time.
Under this situation it is difficult to improve their professional skills
and their teaching levels.

These problems, facing the serious challenges of the 21st century,
cause the inadaptability of the university mathematics education
to the current education development situation to be more serious.
The results are mainly reflected in the following respects.

{\bf For students}, too specific majors and their one-sided understanding of the role of basic courses lead to
their narrow scopes of knowledge (especially knowledge on mathematics), narrow field of vision, lack of creation
and all this is usually expressed in insufficient ``aftereffects''.

{\bf For teachers}, a single- pattern course content,
rigid teaching plans and programs lead to the fact that
they could only be responsible for textbooks and exams and
it is difficult for them to attend to the cultivation of abilities
and qualities of students.
The development of teachers' active roles and their growth are both influenced.

{\bf For textbooks}, most of the content is rather old and their system is stereotyped. There is no mechanism to
encourage teachers to edit new teaching textbooks and materials, which gives people an impression of ``one
thousand people having one face'' and results in difficulties for students to grasp mathematical thinking and
method and to learn new knowledge in mathematics.

\section{Several aspects that should be grasped in doing a good job
of our college mathematics education}

\vskip-5mm \hspace{5mm}

In recent years, with the rapid development of our economy, higher education in our country has met with an
excellent situation of brisk discussions and great reforms. This time of educational reform, whether from the
emphasis degree of the government, or from the realistic spirit of policies, the strength of the funds input and
the depth and breadth of studying problems, they are all unprecedented. Just as the case in previous educational
reforms, the reform of mathematics teaching is still among one of the key problems. What is different from the
past is that the reform is not restricted to the long confusing problem of the relationship between mathematics
and practice, but rather, the mathematics education is first put under the background of social development and is
raised at the height of university quality education. The past is reviewed, the present is examined and the future
is planned. In particular, the role of mathematics education in higher education has been clarified. There emerge
schemes, new textbooks and new courses on mathematics course content, structure and system reforms that have
creative intention and face the epoch. The research and practice on reforming teaching methods of mathematics and
introducing new means of teaching have attracted the attention of broad masses of mathematics teachers; but the
development is not balanced. There is not enough systematic research and comprehensive experiments. Especially,
viewing from the great part of the circle, there are no great changes of the traditional mathematics teaching
modes and methods. As for the renewal of the teaching content of mathematics, it is like ``walking with the
unsteady steps.'' The main problem still lies in the change of educational ideology and teaching concepts. For the
teaching of the basic mathematical courses in non-mathematics class specialties, in our opinion, it is necessary
to clarify the following points in understanding.

{\bf (1)\, It is necessary to have a rather comprehensive
knowledge of the role of mathematics education in university}

In the educational mode in the planned economy, training of professionals in  universities was of an upside-down
type. The national plan decided the arrangement of specialties, and teaching programs were drawn up according to
the requirements of specialties. Courses were arranged according to the procedure of specialty courses -
fundamental specialty courses - basic courses. It was emphasized that the latter shall serve the former. For the
basic education of mathematics, the excessive emphasis of the aspect of ``serving specialties'', and the neglect
of the inherent unity of mathematics as a rational reasoning system and the specific role of mathematics in the
comprehensive quality of students, led to the lack of a comprehensive knowledge of the role of mathematics
education in university education. In fact, mathematics is the common foundation for cultivating and training
various levels of special professionals. For students of non-mathematics class specialties, the role of university
basic mathematics courses lies at least in the following three aspects.

{\bf It is the main course for students to grasp mathematical tools.} Such a role is very important for students
of non-mathematics class specialties and is an important content of ``specialty quality.'' The current problems
are as follows: on the one hand, teachers should study how to effectively enable students to grasp and use this
tool in the whole course of mathematical teaching and how to lay a foundation in this field in the stage of basic
courses; on the other hand, it is necessary to prevent the narrow understanding of ``tool'', viewing the basic
mathematics course only as the tool for serving certain specialty courses, and even the tool specifically for
dealing with certain examinations.

{\bf It is an important carrier for students to cultivate their rational thinking.} What mathematics studies is
the model structures of ``numerals'' and ``forms'', and what it uses is such rational thinking methods as logic,
reasoning and deduction. Large amounts of facts show that it is not a useless theory that is ``isolated from
practice'', but is a thinking creation, which originates from practice and guides practice. The role of such a
training of rational thinking could hardly be replaced by other courses. Such a cultivation of rational thinking
is of utmost importance for students to improve their quality, to enhance their analytical abilities and to
enlighten their creation consciousness. The current problems lie in that we lack the mature teaching experience
and excellent textbooks in this respect.

{\bf It is a way for students to receive the nurture of beauty sense. }
The mathematical aesthetics is part of people's quality of
the appreciation of beauty.
With the development of human civilization and the progress of science,
such a fact is being gradually recognized by people.
In fact, targets mathematics is striving for ---
the arrangement of chaos into order,
the sublimation of experience to laws and the search for the concise and
uniform mathematical expression of motions of substances ---
are the embodiment of mathematical beauty,
and are also man's pursuit for beauty sense.
Such a pursuit acts as a subtle influence for the nurture of
man's mental world and is always a motive force of innovation.
At present we can only say that it is necessary to pay attention
to the role of mathematics in aesthetic education.
And further search and attempts should be made in its embodiment
in teaching and textbooks.

The roles of the three aspects are unified, but in concrete requirements,
according to the different kinds of colleges, specialties and students,
there should be different sides to be emphasized,
so that the comprehensive realization of the knowledge,
abilities and quality in mathematics teaching can be expected.

{\bf (2)\, Emphasis should be put on the solution of problems of
the renewal of college mathematics course system and content }

From the abolition of the imperial civil examination system at the end of Qing Dynasty, the establishment of
modern school and up to the liberation of China, the education of science and technology in colleges and
universities basically followed modes in Europe and America. For the content of mathematical teaching, before the
1950s, mathematics as a compulsory course in most specialties of science and technology except in few specialties,
such as physics, was only limited to simple calculus and differential equations. In the 1950s, the Soviet modes
were comprehensively introduced in higher education institutes. A kind of specialty education with the aim of
training professional talents according to requirements of the trades was formed. The arrangement of the courses
in the specialties was aimed at the requirements of the professional knowledge. Basic courses such as mathematics
were required to serve specialty courses. Such a system and guiding ideas have actually continued till now. Though
in the 1980s, due to the need to employ computers, linear algebra (taking calculation as the chief content) and
numerical methods became compulsory mathematics courses, viewing from the guiding ideas, the main point was still
``specialty education'' and the changes of teaching content were very limited. Later, there occurred improper
evaluations, in addition to the influences of ''examination-oriented education'' induced by taking the
post-graduate entrance examinations, the tendency of one-sided pursuit of the problem-solving techniques was
promoted. As a result, the train of thought about reforms became more vague. However, due to the rapid development
of the computer technique in the later half of last century, people had more knowledge on the role of modern
mathematics. Especially in the recent 30 years, there appears the so-called ``modern mathematics technique'' which
displays its prowess fully in economy and industries. For example, optimization, engineering control, information
processing, fuzzy recognition and image reconstruction, etc. they are produced due to the combination of the
principles and methods of modern mathematics with computers. They penetrate and are applied in all sectors and
trades, and are combined with relevant techniques to form the so-called high and new techniques in those areas. In
current developed countries, mathematics is applied to improve the organizational level of economy, and from
drawing up macroscopic strategic planning to the storage, distribution and transportation of products and to the
market prediction, analysis of finance and insurance businesses, significant progress has been attained. All the
facts mentioned above mean that mathematics has turned from part of the traditional natural sciences and
engineering techniques to further penetrating into many areas of the modern society and economy and has gradually
become one of their indispensable columns. On the other hand, the development of science and technology brings
about a series of problems and it is necessary for people to reexamine the relationship between man and society
and between man and Nature in a rational way. All this requires that training on the tool property and rationality
should be emphasized in college mathematics education. Therefore, it is an urgent task to adjust the basic course
system of college mathematics and appropriately renew teaching content.

{\bf (3)\, Reforming the examination-oriented
``instillation-type'' teaching}

In recent years, various kinds of factors have induced a type of examination-oriented teaching. Marks play almost
decisive roles for students in many respects such as evaluation in the prizes, graduate entrance examinations and
the assignment for jobs. In some occasions, students' marks are the main factor for evaluating teachers. Such a
situation causes teachers to teach for exams and students study for exams. In order to raise the average marks of
a class, teachers spend a lot of time and energy on the teaching of problem patterns. Accordingly, the
``instillation-type'' of teaching methods are in a dominating position in our college teaching. Characteristics of
such methods are as follows: on the one hand, various kinds of problem patterns are explained in every detail in
the classroom and efforts are made to enable students to understand them as soon as they listen to the
explanations. In such a way, besides the small amount of classroom information quantity, the psychology of
dependence of the students will certainly be enhanced and they will have a bad habit of being lazy in thinking,
which will seriously hinder the cultivation of the innovation sense and innovation capability. On the other hand,
students are busy problem solving after the class and problems are substituting learning. They rarely read books
and pay little attention to the mathematical thinking and mathematical applications, let alone the cultivation of
the innovation sense.

{\bf (4)\, Devoting great efforts to the construction of
teachers' teams}

This is the crux matter of our reform.
In recent years, compared to the past, the situation of teachers' teams
engaged in non-mathematics class specialty mathematics teaching has turned
to the better. Quite a few of young teachers
who have obtained Doctor's degrees join the career.
In common colleges and universities, proportions of mathematics teachers
who have gained Master's degrees are increasing.
This is an encouraging phenomenon, indicating that there will be qualified
successors to carry on mathematics teaching, but the authorities of colleges
and universities have not timely put such a strategic task at a deserved height.
For example, knowing clearly that the experienced teachers
with high academic levels should be appointed to teach the basic courses,
some colleges and universities still classify the teachers according
to the order of the postgraduate teaching, senior students teaching,
and the basic courses teaching. Such a guiding direction is
very unfavorable for cultivating high-quality talents.
In addition, viewing from the young teachers' teams,
there also exists something worrying. First, their ideas are not so stable;
secondly, they cannot deal well with the relationship
between scientific research and teaching.
In the key universities, there usually exists the atmosphere of
paying more attention to scientific research than that to teaching.
In common colleges and universities, there always exists the problem
of only teaching without scientific research.
Therefore, creating certain conditions and asking the young teachers
to participate in scientific research are important measures
to improve the levels of college mathematics teachers.
Thirdly, not enough efforts have been made in conducting education
on teacher morality and on superior teaching traditions for the young teachers.

\section{Several important relationships
in the university mathematics education reforms}

\vskip-5mm \hspace{5mm}

The college mathematics education reform is a very careful job.
The treatment of many problems in this respect will
not only rely on principles, but also be flexible.
According to the experience and lessons from
the past mathematics education reform,
I would like to discuss some problems concerning principles
in mathematics teaching reform from the correct treatment of
the relationships in several aspects as follows:

{\bf (1)\, Embodying quality education, paying attention to well
treating the relationships between knowledge impartment and
quality cultivation}

In recent years, there are a large amount of discussions on ``quality'' and ``quality education''. We assume that
quality is a mental and physical attribute expressed when people understand and treat things and events. It is
based on the congenital psychological conditions and is gradually formed under the influences of the postnatal
environment. Viewing from the angle of education, an individual's quality in a certain aspect is shown in his
power of understanding and potential energy in such an aspect. Quality education is a process during which
excellent qualities of people are cast through the smelting of systematic knowledge. Any objective existence in
the world has its attribute characters of ``numerals'' and ``forms''. With the progress of the science of
mathematics, people's knowledge of ``numerals'' and ``forms'' has enhanced from directly-visual quantity
relationships and space forms to the abstract ``mathematical structures'' and ``space concepts'' with deeper
connotation and wider extension. Such attributes of ``numerals'' and ``forms'' in people's understanding things
and their comprehension and potential ability in handling the corresponding relationships are obviously a kind of
people's quality and we call such a quality mathematical quality. Under the background of the tendency of
digitalization and informationalization in the knowledge economy society, the importance of possessing such a
quality is very great. In the quality education, viewing from the main body of education (i.e. the main targets of
education and corresponding teaching behaviors), the soul of the college mathematics education is exactly the
mathematical quality, and is what has been discussed in the previous part --- quality property. Viewing from
mathematics itself, especially the modern mathematics, its essence is a rational thinking system extracted
abstractly from large amounts of objective phenomena. Therefore, in its education, besides displaying how it
absorbs nutrition from laws of vivid objective things and provides tools, main focus is the training of the
rational thinking techniques and of the mastering mathematical tools. That is the content of the education of
mathematical quality. Though it imparts ``knowledge'' as well, ``knowledge'' of mathematics acts as a carrier of
imparting quality (such a role was always neglected by people in the past), besides one side that it is combined
with the material content and thus acts as a tool to serve other disciplines. The quality education in the
mathematical teaching is that the teachers put the lively and rational ways of thinking via the knowledge carrier
to implement the motivational, psychological and mental guidance for their students.

Mathematics has become a fundamental component of the culture of the contemporary social culture via its ability
of tools, rational spirit and sense of beauty. In the society of the 21st century if one has no idea what
mathematical technology means, lacks both the rational thinking and the sense of beauty appreciation; then his
total quality will be affected. His abilities in insight, judgement and originality would be greatly restricted.
Therefore, the accomplishment of mathematical culture is not a kind of ``fashion'' for one in the
ever-increasingly fierce competition among talents of the society; it is indeed a kind of actual necessity for one
in his work, study and social communication.

The specific content of the mathematical quality is quite rich,
the most outstanding specific characteristics that are
generally acknowledged are summarized as follows:

\begin{itemize}
\item[$\bullet$]
The sharp consciousness of extracting the attributes of ``numerals and forms'' of things.
\item[$\bullet$]
Thinking mode of using ``abstract models and structures'' to study things.
\item[$\bullet$]
Exploring habits of conducting compact deduction
by means of signs and logical systems.
\end{itemize}

Features of mathematical qualities in these aspects are connected with
each other and cannot be separated in the actual educational processes.
All college students should be cultivated and educated in such aspects.
Of course, the education and cultivation of mathematical qualities
are different in emphasis points, depth, and breadth between
mathematics major students and non-mathematics class specialties students,
and in non-mathematics class specialties, such differences do exist
between students of science, engineering, agriculture and medicine
and those of humanities majors.

{\bf (2)\, Grasping course systems and content renewal and
treating well relationships between the inherence and
modernization of the mathematical knowledge}

The main body of the textbooks of current basic mathematics courses is mostly mathematics before the 19th century.
This is in contrast with other basic textbooks of physics, chemistry and biology that began mostly from the 19th
century. Therefore, the problem of the modernization of basic university mathematics textbooks is to be stressed.
But we should be cautious in this problem. It is unavoidable to appropriately add the commonly accepted and
fundamental content of modern mathematics; but it is necessary to consider another aspect: as we have mentioned
above, mathematics is a kind of thinking science and has its own specialties. Its system is constructed by logic
and is a structure of series established one layer after another from the bottom to the top. The mathematical
knowledge important in the past is the ``logical basis'' of the current mathematics and throwing away the former
will influence the later studies. For example, calculus has a history of over three hundred years, but it still is
the foundation stone of modern mathematics and can't be thrown away casually. Such a situation differs from other
sciences. For other sciences, if the theories before the 19th century are out-of-date, they can be cancelled,
though proceeding and subsequent knowledge has inherent property. The preceding knowledge is not necessarily the
direct logical foundation of the subsequent knowledge, so the abandonment of the preceding knowledge will not
influence the learning of the subsequent knowledge. In addition, there is still an important reason. Such a part
of mathematics content, for example, calculus, has still a rather wide application today. Therefore the idea of
``modernization'' of the content of college mathematics textbooks is somewhat different from the modernization of
other disciplines and thus how much modern content has been listed in the textbook cannot be simply taken as the
standard of measurement. We assume that the content modernization of the college mathematics textbooks can mainly
include the following aspects. First, the classical mathematics content should be governed by viewpoints and
languages of modern mathematics as far as possible; at the same time ``roots'' in the classical mathematics for
some modern mathematics should be appropriately introduced. Secondly, significant results of modern mathematics
that have formed the basic part of relevant disciplines should be put into textbooks as far as possible and a
poplar introduction should be made. Lastly, it is necessary for students to possess essential basis of modern
mathematics necessary for students to further study relevant specialties by themselves.

What is similarly important is to cancel some contents such as reasoning loaded down with trivial details and
those calculations that can be done by calculators, and concepts and methods that are relatively old and thus have
no development prospect in modern science. In summary, facing an accumulation of nearly two centuries, how to deal
with the ``metabolism'' of the basic content of mathematics is a quite difficult problem. Depending only on
``outside extension type'' to increase academic hours will not work and efforts should be made on reforms on the
structure and content of the course.

{\bf (3)\, Paying attention to the reform of teaching methods and
handling well relationships between teachers' guidance and
students' initiative}

Once the teaching system and content of the course of college mathematics
have been determined, teaching methods become the key problem
of teaching quality. Obviously, the teaching process of basic mathematics
is absolutely not man's repetition of cognition processes for numerals
and forms laws, but basic laws of cognition processes should be observed.
In addition, the teaching process cannot be the simulation
of research process of mathematical problems,
but should embody research methods and thinking modes specific for mathematics.
The teaching of mathematics includes two sides of teaching and
learning and is a comprehensive process in which many teaching activities
guided by teachers lead to the initiative of students to study mathematics.
The so-called initiative has at least the following meanings.
First, students have the interest and dynamic force to learn mathematics.
Secondly, students can learn by themselves under the proper guidance of teachers.
Thirdly, students can work at and put forward problems independently
during the study process. Fourthly, students can use
the knowledge their teachers impart and related knowledge gained
through self-study after processing and digestion to construct
a knowledge system that, even if it is simple and not perfect,
actually belongs to themselves.

{\bf (4)\, Emphasizing the practical links of mathematics and
paying attention to handling well relationships between
mathematical basic training and mathematical application
consciousness and ability cultivation}

The widespread use of computers and the emergence of a series of powerful mathematical software systems have
resulted in profound changes in the roles of mathematics and in mathematical teaching. They make it possible to
collect and process large amounts of data, also make mathematical model a means of experiment, and thus greatly
promote applications of mathematics in all fields. The combination of mathematical thinking with computers has
become an important mode of modern mathematical teaching. For example, since the establishment of the course of
``mathematical models'' in colleges in the 1980s, there have been hundreds of colleges and universities where this
course has been set up, and the course is a favorite of students. Especially the annual national contest of
college mathematics model-building attracts students of various specialties and promotes reforms in mathematics
teaching. At present, the course ``mathematical experiments'' aiming at strengthening practice links of
mathematics is being tested in several colleges and universities and the preliminary results are encouraging.

\label{lastpage}

\end{document}